\documentclass{amsart}
\usepackage{amsthm,amssymb}
\newtheorem{theorem}{Theorem}

\newtheorem{lemma}[theorem]{Lemma}

\newtheorem{remark}{Remark}

\newcommand{\CU}{\mathcal{U}}
\newcommand{\TI}{\mathnormal{T} \mathnormal{I}}

\newcommand{\CP}{\mathcal{P}}
\newcommand{\CQ}{\Delta \mathcal{P}}

\begin{document}
\title[Tur\'{a}n inequalities ]{Tur\'{a}n inequalities and  zeros of orthogonal polynomial}
    
\author[I. Krasikov]{Ilia Krasikov}
     
\address{  Department of Mathematical Sciences,
	   Brunel University,
	   Uxbridge UB8 3PH, United Kingdom}
\email{mastiik@brunel.ac.uk}
			      
			       
%\date{}

 %-------------------- Title Page ----------------------------------
  
\begin{abstract}
We use Tur\'{a}n type inequalities to give new non-asymptotic bounds on the extreme
zeros of orthogonal polynomials in terms of the
coefficients of their three term recurrence.
Most of our results deal with symmetric polynomials satisfying
the three term recurrence $p_{k+1}=x p_k-c_k p_{k-1},$
with a nondecreasing sequence $\{c_k\}$. As a special case
they include a non-asymptotic version of M\'{a}t\'{e}, Nevai and Totik
result on the largest zeros of orthogonal polynomials with
$c_k=k^{\delta} (1+ o(k^{-2/3})).$

\end{abstract}
   
   %------------------------------------------------------------------------%
   %                                                                        %
   %         1. INTRODUCTION                                                %
   %                                                                        %
   %------------------------------------------------------------------------%
\maketitle
\section{Introduction}

Let $\CP = \{ p_i(x) \} $ be a family of orthogonal polynomials 
satisfying the three term recurrence relation
\begin{equation}
\label{rec}
b_k p_{k+1}(x)=(x-a_k)p_k(x)-c_k p_{k-1}(x), \; \; \; b_k, c_k >0,
\end{equation}
with the initial conditions $p_{-1}=0,\; p_0=1,$
and let $x_{1k} <...< x_{kk} ,$
be the zeros of $p_k(x).$ 
We are interested in finding uniform bounds on the extreme zeros, that is
an interval $I=[A(k),B(k)],$ such that $A(k) <x_{1k}<x_{kk}<B(k),$
in terms of the coefficients of the recurrence.
Such a setting arises naturally if one deals with the family 
depending on parameters, as in the case of classical Jacobi and
Laguerre polynomials, and is seeking for bounds uniform in all the parameters involved. 
The main aim of this paper is to show that the classical Tur\'{a}n inequality
\begin{equation}
\label{tur}
T_k (\CP,x)=p_k^2(x)-p_{k-1}(x)p_{k+1}(x) \ge 0,
\end{equation}
and its analogues (abbreviated as $\TI$ in the sequel),
provide a convenient tool for tackling the problem.
It is known that (\ref{tur}) holds for some families of
orthogonal polynomials including Laguerre, Jacobi and some other polynomials
(see \cite{szwarc} and the references therein).
\\
At present there are two general approaches to the problem, 
one is based on the chain sequences \cite{ismail} and another exploiting
the Rayleigh quotient to find the extreme eigenvalues
of the corresponding Jacobi matrix. The last one was discovered independently
by G. Freud \cite{freud} and V.I. Levenshtein 
(see \cite{leven}, where the references on the earlier papers of the author are
given), and yields the following elegant representation for the extreme zeros.
\begin{equation}
\label{levshm}
x_{1k}=\min \left(\sum_{i=0}^{k-1}a_i x_i^2-2 \sum_{i=0}^{k-2}x_i x_{i+1} \sqrt{b_i c_{i+1}} \right),
\end{equation}
\begin{equation}
\label{levshM}
x_{kk}=\max \left(\sum_{i=0}^{k-1}a_i x_i^2+2 \sum_{i=0}^{k-2}x_i x_{i+1} \sqrt{b_i c_{i+1}} \right),
\end{equation}
where the extrema are taken over all (or only over positive) $x_0,x_1,...,x_{k-1},$ subjected to
$\sum_{i=0}^{k-1} x_i^2 =1.$
\\
For symmetric polynomials (i.e. when $p_k (-x)=(-1)^k p_k(x)$)
and the monic normalization, the case 
we mainly deal with in this paper, the recurrence
(\ref{rec}) can be rewritten as
\begin{equation}
\label{recsym}
p_{k+1}(x)=x p_k(x)-c_k p_{k-1}(x),
\end{equation}
and (\ref{levshm}), (\ref{levshM}) become
\begin{equation}
\label{levglav}
x_{kk}=-x_{1k}=2 \; \max \left( \sum_{i=0}^{k-2}x_i x_{i+1} \sqrt{c_{i+1}} \right).
\end{equation}
Thus $x_{kk}$ as a function of the vector $(c_1,c_2,...),$ possesses many nice properties, e.g.
it is clearly subadditive, continuous, increasing with $k.$ This makes transparent many otherwise puzzling
questions concerning the behaviour of the extreme zeros.
The most striking result
obtained via (\ref{levshm}), (\ref{levshM}) is due to
A. M\'{a}t\'{e}, P. Nevai and V. Totik \cite{mtn} and states
that if $c, \delta >0,$  are fixed and
$c_k=c^2 k^{2 \delta} \left( 1+o(k^{-2/3}) \right), $ 
then
\begin{equation}
\label{nevoz}
x_{kk} k^{-\delta} = 2 c -c \cdot 3^{-1/3}(2 \delta )^{2/3} i_1 k^{-2/3}+o(k^{-2/3}),
\end{equation}
where $i_1=3.3721...,$ is the smallest zero of the Airy function.
\\
If $x_{kk}$ is represented by a series in decreasing powers of $k,$
it is naturally to distinguish between the first and the second order bounds,
e.g. $2c$ and $2c(1-K \delta ^{2/3} k^{-2/3})$ in (\ref{nevoz}),
where to maintain the uniformity in $k$ 
we will allow a weaker constant $K$ than the exact asymptotic one.
Whenever the first order bounds can be obtained rather easily, say,
by replacing $2 x_i x_{i+1}$ in (\ref{levglav}) by $x_i^2+x_{i+1}^2 ,$ or by similar
elementary arguments \cite{assche00},\cite{leven},
the second order estimates uniform in all the parameters, which are the main subject of this
paper, were found only recently
for the case of classical orthogonal polynomials \cite{klag} - \cite{kf}.
\\
It is worth also noticing that 
it is rather easy to extract bounds on the zero
$x^*_{kk}$ corresponding to  
the perturbed recurrence $p_{k+1}(x)=x p_k(x)-c_k (1+ \epsilon _k)^2 p_{k-1}(x),$
provided one knows the result for $\epsilon =0.$ 
Indeed, 
$\sum_{i=0}^{k-2} x_i x_{i+1} \le 1,$
and thus, 
$x^*_{kk}$ is in the interval $(1 \pm \epsilon ) \; x_{kk},$
for $| \epsilon_k | <\epsilon < 1 .$
For example it would be enough to establish (\ref{nevoz}) for $c_k=c^2 k^{2 \delta},$
the general case with the extra factor $( 1+o(k^{-2/3})),$ follows from the above arguments.
In this paper we will not state explicitly such obvious generalizations, 
but the reader should keep them in mind. 
\\
Our approach to the problem is quite different from the above two
and based on the following observation \cite{k3}.
Let $f=f(x)$ and $g=g(x),$ $ |\deg (f)-\deg (g)| \le 1,$ be two real polynomials 
with only real interlacing zeros. Suppose also that we know a
(preferably quadratic) form 
$ \sum_{i=0}^m A_i(x)f^{m-i}g^i \ge 0,$ which is indefinite in $f$ and $g$ viewed as
formal variables. 
Then one can routinely obtain
bounds 
on the extreme zeros of $f$ and $g $ similar to (\ref{nevoz})
but, of course, with a weaker constant instead of $i_1.$
The existence of such a form is far from being obvious
and in fact the main difficulty is to find an appropriated one.
For classical orthogonal polynomials,
when a second order differential equation is known,
one may choose $f=p_k, \; g=p'_k,$ and the quadratic form obtained
from the Laguerre inequality ${f'}^2-f f'' >0,$ or its higher order generalizations
\cite{k3, klag, kest, kf}. In the discrete case $f=p_k(x+1), g=p_k(x),$ play a similar role
\cite{k3,kd,kch}.
\\
In this paper we will show that the required forms, in particular these giving second
order bounds, may be
obtained directly from  $\TI$
via the three term recurrence.
To this end we will establish two new sets of $\TI$ 
yielding second order bounds.
The following theorem is one of our main results.
\begin{theorem}
\label{th1}
Let $p_k$ be a symmetric polynomial
satisfying (\ref{recsym}) and suppose that $c_k$ are nondecreasing, then
$$x_{kk} < 2 \sqrt{c_{k-1}}, \; \;  k \ge 2. $$
Moreover if  $d_i =\frac{c_i -c_{i-1}}{c_i} \ge 0, \; \; c_0=0,$
satisfy
\begin{equation}
\label{condnew112}
\frac{d_i}{2(1+d_i)} < d_{i+1} < \frac{d_i(1+2 \sqrt{d_i}+2 d_i )}{1+d_i}\; ,
 \; \; \; i=1,2,...;
\end{equation}
then for $k \ge 2,$ 
\begin{equation}
\label{condsimpl}
x_{kk} < 2  \;  \sqrt{ c_k \left( 1- \frac{d_{k+1}^{2/3}}{(2^{1/3}+d_{k+1}^{1/3} \; )^2} \right) } 
\; , \; \; \; k \ge 2.
\end{equation}
\end{theorem}
As a corollary we obtain the following uniform version of (\ref{nevoz}).
\begin{theorem}
\label{thmnt}
Let $c_k=c^2 k^{2 \delta} ,$ where $c >0, \; \delta  \ge 0,$ are fixed,
then
$$x_{kk} k^{-  \delta} < 2c \; \sqrt{1-\frac{\delta^{2/3}}{((k+
\frac{1}{2})^{1/3}+\delta^{1/3} )^2} } \;, \; \; \; k \ge 2 .
$$
\end{theorem}
It would be important to have lower bounds corresponding to (\ref{condsimpl}).
A trivial one $x_{kk} > \sqrt{c_{k-1}}$ readily follows from (\ref{levglav}) (see also \cite{assche00},
\cite{kf}) in this connection).
One hardly could expect that (\ref{condsimpl}) is sharp for a
faster than polynomial
rate of growth of $c_k,$ 
although formally Theorem \ref{th1}  allows $c_k$ to be of order
about $e^{k^2/2}.$ It would be also interesting to obtain similar results for
polynomials orthogonal on $[-1,1]$ and $[0, \infty ).$

The paper is organized as follows.
In the next section we survey some known $\TI$ which will be used
in section 3 for obtaining bounds on the extreme zeros.
In the last section we establish two new sets of $\TI$ and
give second order bounds for a vast class of orthogonal polynomials satisfying (\ref{recsym})
with a nondecreasing sequence $c_k.$

\section{Tur\'{a}n Inequalities}
The reason for study $\TI$
in
the theory of orthogonal polynomials
is that they have
a few quite important applications.
For example,
under appropriate restrictions
$T_k ( \CP,x) $ converges uniformly on compact subsets of $(-1,1)$
to $\frac{2 \sqrt{1-x^2}}{\pi \alpha'(x)},$
where $\alpha'(x)$ is the absolutely continuous part of the corresponding orthogonality
measure on $[-1,1].$ The existing technique enables one to obtain similar limiting expressions
for higher order analogues of $\TI$ as well, e.g.
the right hand side of inequality (\ref{marik}) below tends to
$\frac{8(1-x^2)}{\pi^2 (\alpha'(x))^2},$
see \cite{dimit} and the references therein.
\\
Today there are two almost independent theories related to $\TI$,
the first dealing with the validity of (\ref{tur}) for different families
of orthogonal polynomials and goes back to
the pioneer work of Tur\'{a}n \cite{turan}.
Another one is motivated by the theory of entire functions, in particular by the Riemann hypothesis.
In the last case much more precise higher order 
generalization of (\ref{tur}) are known, but 
the orthogonal polynomials must have a generating function of a very
special type.
A survey of this theory and the relevant references
can be found in \cite{craven},\cite{craven1},\cite{dimit}.
\\
Progress in the first direction was recently summarized by R. Szwarc \cite{szwarc},
who established the following result.
\begin{theorem}
\label{thswar}
\noindent
$(i)$ Let $\CP$ be a family of symmetric
polynomials orthogonal on $[-1,1],$
where three term recurrence (\ref{rec}) is normalized by
\begin{equation}
\label{symfin}
b_k+c_k=1, \; c_0=0, \; c_{k+1} >0, \; b_k >0 ,
\end{equation}
that is by $p_k(1)=1, \; p_k(-1)=(-1)^k .$
Then
\begin{equation}
\label{tur11}
T_k (\CP,x) \ge 0, \; \; |x| \le 1,
\end{equation}
provided one of the following conditions holds
\\
$(i_a)$  $c_k$ is nondecreasing and $c_k  \le \frac{1}{2};$
\\
$(i_b)$ $c_k$ is nonincreasing and $c_k  \ge \frac{1}{2}.$
\\
$(ii)$ Let $\CP$ be a family of
polynomials orthogonal on $[0,\infty),$ which are normalized by $p_k(0)=1,$
with the three term recurrence 
\begin{equation}
\label{halfinf}
x p_k(x)=- b_k p_{k+1}(x)+(b_k+ c_k ) p_k(x)-c_k p_{k-1}(x) ,
\end{equation}
where $b_0=1, c_0 =0, b_k >0, \; c_{k+1}>0.$
Suppose $b_k$ and $c_k$ are nondecreasing, then
\begin{equation}
\label{tur12}
T_k (\CP,x) \ge 0, \; \; x \ge 0,
\end{equation}
provided one of the following conditions holds
\\
$(ii_a)$ $c_k \le b_k , \; \; c_k - c_{k-1} \ge b_k- b_{k-1};$
\\
$(ii_b)$ $c_k \ge b_k , \; \; c_k - c_{k-1} \le b_k- b_{k-1}.$
\end{theorem}
R. Szwarc \cite{szwarc} also obtained similar yet rather technical
conditions which guarantee the validity of (\ref{tur}) for a general 
nonsymmetric polynomial with  the orthogonality measure supported on $[-1,1].$
For the sake of simplicity we 
did not state them here.
It seems nothing is known in the case of  nonsymmetric polynomials orthogonal
on the whole real axis.
On the other hand
the case of symmetric polynomials orthogonal on $(- \infty , \infty )$ is almost trivial.
\begin{theorem}
\label{trivial}
Let $\CP$ be a family of orthogonal
polynomials satisfying (\ref{recsym}) with a nondecreasing sequence $\{ c_k \}_{k=1}^\infty .$
Then $T_k (\CP,x) \ge 0.$
\end{theorem}
\begin{proof}
The result follows by
$T_0 =1,$ and an easy to check identity
\begin{equation}
\label{firsto}
T_{k+1}(\CP,x)=c_k T_k (\CP,x)+(c_{k+1}-c_k)p_k^2(x).
\end{equation}
\end{proof}
A few higher order generalization of (\ref{tur}) are known. 
To state them we 
shall consider the Laguerre- P\'{o}lya class of functions which consists of
real polynomials with only real zeros and real entire functions
$$F(z)=c e^{-\alpha z^2+ \beta z } z^r \prod_i (1-z/z_i)e^{z/z_i},$$
where $\alpha \ge 0, \; c, \beta, $ are real, $r$ is a nonnegative integer and
$\sum_i z_i^{-2}$ is convergent. 
\\
Suppose now that a family $\CU$ of real functions
$u_k=u_k(x), \; k=0,1,...,$  has for some values of $x$
a generating function of the Laguerre- P\'{o}lya class, 
\begin{equation}
\label{pl}
\sum_{k=0}^{\infty} u_k \frac{z^k}{k!} =F(z).
\end{equation}
An instructive example is provided by
the binary Krawtchouk polynomials $u_k=K_k^n(x),$ $n$ is a positive integer,
having the generating function
$$\sum_{i=0}^\infty K_i^n(x) z^i =(1-z)^x (1+z)^{n-x},$$
which satisfy (\ref{pl}) for $x=0,1,...,n.$
Another examples are given by ultraspherical $C_k^{( \lambda )},$ Hermite $H_k ,$
and Laguerre $L_k^{( \alpha )} $ polynomials \cite{dimit},
with
$$u_k=\frac{C_k^{( \lambda )}(x)}{C_k^{( \lambda )}(1)}, \; \; \lambda >-\frac{1}{2}, \; \; -1 \le x \le 1;$$
$$u_k(x)=H_k (x), \; \; -\infty <x< \infty ;$$
$$u_k=\frac{L_k^{( \alpha )}(x)}{L_k^{( \alpha )}(0)}, \; \; \alpha >-1, \; \; x \ge 0.$$
\\
In the following
theorem the first part belongs to M. Patrick \cite{patrick2} 
and the second one to J. Ma\v{r}\'{\i}c \cite{marik} 
(the extension of it to the whole
Laguerre- P\'{o}lya class is due to D.K. Dimitrov \cite{dimit}).
\begin{theorem}
For those values of $x$  for which (\ref{pl}) holds
\begin{equation}
\label{patr}
T^{(m)}_k(\CU ,x) = \frac{1}{2} \; \sum_{j=0}^{2m} (-1)^{j+m} {2m \choose j} u_{k-m+j}u_{k+m-j} \ge 0, \; \;
m=0,1,... \; ,
\end{equation}
and
\begin{equation}
\label{marik}
S_k(\CU,x) = 4(u_k^2-u_{k-1}u_{k+1})(u_{k+1}^2-u_k u_{k+2}) - (u_k u_{k+1}-u_{k-1}u_{k+2})^2  \ge 0.
\end{equation}
\end{theorem}
Notice that 
$T^{(1)}_k(\CU,x)$  
is just $T_k (\CU,x).$
The inequality (\ref{marik}) can be viewed as a refinement of 
(\ref{tur}) and is intimately connected with so-called Newton inequalities \cite{nicu}.
To apply 
(\ref{patr}),(\ref{marik}) to orthogonal polynomials it would be important to
to restate the condition for their validity
in terms of the coefficients of (\ref{rec}).
For $T^{(2)}_k$ and recurrence (\ref{recsym}) this will be accomplished 
in the last section. The corresponding question for
$S_k$ remains open.
%%%%%%%%%%%%%%%%%%%%%%%%%%%%%%%%%%%%%%%%%%%%%%%%%%%%%%%%%%%%%%%%%%%%%%%%%%%%%5555555
\section{Extreme zeros}
In this section we give some first order bounds on the extreme zeros
which can be deduced from
Theorems \ref{thswar} and \ref{trivial}.
We also show that $S_k(\CU,x)$ (which gives a fourth degree form),
yields second  order bounds for Hermite polynomials.
The inequality $T^{(2)}_k \ge 0,$ will be considered in more details
in the next section. 
\\
First, we describe simple geometric arguments which enable one to
deduce bounds on the extreme zeros from a given form.
Let $p=p(x), \; q=q(x)$ be two real polynomials, 
$\deg(p)=k \ge 2, \; \deg(q)=k-1,$
with only real interlacing zeros
$x_1,...,x_k ,$ and $y_1, ...,y_{k-1},$  respectively, $x_i<y_i <x_{i+1}.$
Suppose that for $x \in (M,N), \; M<x_1 <x_k <N,$ where $M$ and $N$ can be finite or infinite,
there exists a nonnegative form
\begin{equation}
\label{kvfr0}
\sum_{i=0}^m A_i(x) q^i p^{m-i} \ge 0, \; \; m \ge 2, \; \; even,
\end{equation}
where $A_i (x)$ are certain functions defined in all the points of $(M,N).$
Introducing  the function $t=t(x)=q/p,$ we rewrite it as
$$
Q(t,x)=\sum_{i=0}^m A_i(x)t^i \ge 0.
$$
Since
$\lim_{x \rightarrow \pm \infty} t(x)=0,$
and the zeros of $p$ and $q$ are interlacing
then
$t(x)$
consists of two hyperbolic $B_0,B_k$  and $k-1$ cotangent-shaped decreasing branches $B_1,...,B_{k-1},$
where $B_i$ is defined for $x_{i}<x<x_{i+1},$ $x_{0}=-\infty, \; x_{k+1}=\infty .$
A function $t_0=t_0(x)$ will be called an $(i,j)-transversal$ if 
\\
(i) $t_0$ is continuous on $[M,N],$
\\
(ii) $t_0$ intersects each of the branches $B_l, $ of $t$ for $ i\le l \le j.$ 
\\
(iii) there is an open interval $I \subset [M,N]$ such that $Q(t_0,x) \le 0$ iff $x \in [M,N] \backslash I.$

Obviously, if an $(i,j)-$transversal exists then $[x_{i+1}, x_j] \subset I,$ 
and to get bounds on $x_{i+1},x_j$ one needs just to find the extreme roots of the equation
$Q(t_0,x)=0,$ on $[M,N].$ 
Note that any continuous function intersects all the cotangent-shaped branches
$B_1,...,B_{k-1},$
and is, if (iii) holds, a $(1,k-1)-$transversal, thus giving bounds on $x_2$ and $x_{k-1}.$
For example, $t_0=c x,$ is 
a $(0,k)-$transversal for $c >0,$ and 
a $(1,k-1)-$transversal for $c \le 0,$ provided it satisfies (iii).
As we will see, in many cases the condition (iii) is automatically fulfilled and moreover, a naive 
choice of $t_0$ 
as a solution of $\frac{\partial Q(t,x)}{\partial t}  =0,$
that is as the function providing the minimum to $Q(t,x),$ does work.
The situation is  especially simple for quadratic forms, the case we mainly exploit here.
Then $Q(t,x)=A_0+A_1 t+A_2 t^2,$ and one may try $t_0= -\frac{A_1}{2A_2},$
with $Q(t_0,x)=\frac{4A_0 A_2 -A_1^2}{4A_2}.$
\\
In the rest of the paper we will use $t=p_{k-1}/p_k ,$
and with one explicitly stated exception,
$t_0$ will be chosen as a solution of $\frac{\partial Q(t,x)}{\partial t}  =0.$ 
\\
The simplest way to obtain the required quadratic form $Q(t,x)$
for orthogonal polynomials
is to   
express $p_{i-1}, p_i,$ and $ p_{i+1}$ in $T_{i}, \; |k-i| \le 1,$ via 
$p_{k-1},p_k$ by
the three term recurrence.  
In this case
one gets three (slightly different) bounds on the zeros, we present just one of them in the
theorem below.
But already for $T_{k \pm 2},$
the above expression for $t_0$ may have singularities and our arguments
are not applicable without certain  restrictions on the coefficients
(see Lemma \ref{tt2} below).
Using Theorems \ref{thswar} and \ref{trivial}
to  guarantee the corresponding $\TI$ 
we get the following first order bounds.
\begin{theorem}
\label{firstob}

\noindent
$(i)$ Let $p_k$ be a symmetric polynomial
satisfying (\ref{recsym}) and suppose that $c_k$ are nondecreasing, then
$$x_{kk} < 2 \sqrt{c_{k-1}}, \; \;  k \ge 2. $$
\noindent
$(ii)$ Let $p_k$ be a symmetric polynomial orthogonal on $[-1,1]$ satisfying 
(\ref{rec}) and (\ref{symfin}),
then
$$
|x_{ik}| < 2 \sqrt{ b_k c_k },$$
where $i=1,...,k,$
if $c_k$ is nondecreasing and $c_k  \le \frac{1}{2};$
and $i=2,...,k-1,$ if $c_k$ is nonincreasing and $c_k  \ge \frac{1}{2}.$

\noindent
$(iii)$ Let $p_k$ be a polynomial orthogonal on $[0,\infty)$ satisfying (\ref{halfinf}).
If $b_k$ and $c_k$ are nondecreasing, then
$$(\sqrt{b_k }-\sqrt{ c_k} \; )^2<x_{2,k} < x_{k,k}<( \sqrt{b_k }+\sqrt{ c_k} \; )^2,$$
provided $c_k \le b_k , \; \; c_k - c_{k-1} \ge b_k- b_{k-1};$
and
$$(\sqrt{b_k }-\sqrt{ c_k} \; )^2<x_{1,k} < x_{k,k}<( \sqrt{b_k }+\sqrt{ c_k} \; )^2,$$
provided  $c_k \ge b_k , \; \; c_k - c_{k-1} \le b_k- b_{k-1}.$
\end{theorem}
\begin{proof}
$(i)$ By (\ref{recsym}) we get
$$Q(t,x)=c_{k-1}p_k^{-2} T_{k-1}(\CP,x) =1-x t +c_{k-1}t^2  \ge 0.$$
In our case $M=N= \infty,$ 
and $t_0= \frac{x}{2c_{k-1}},$ is clearly a $(0,k)-$transversal.
Finally $4c_{k-1}^2 Q(t_0,x)=4c_{k-1}-x^2,$ and the result follows.
\\
$(ii)$ Substituting $p_{k+1}$ from (\ref{rec}) (in our case $a_k=0$), we have for $[M,N]=[-1,1],$
$$Q(t,x)=b_k T_k (\CP,x)=b_k -x t+c_k t^2 \ge 0, $$ 
with
$t_0=\frac{x}{2c_k},$
and $4b_k c_k  Q(t_0,x)=4b_k c_k-x^2.$
Obviously, $t_0$ is
a $(1,k-1)-$transversal.
Finally, using the normalization $p_i(1)=1, \; p_i(-1)=(-1)^i,$ and $b_k+c_k=1,$
one can check that $t_0(-1) \le t(-1),$ and $t_0(1) \ge t(1),$ only if $c_k \le \frac{1}{2}.$
Thus, for $c_k \le \frac{1}{2},$ $t_0$ intersects all the branches of $t$ and hence is 
a $(0,k)-$transversal.
\\
$(iii)$ Substituting $p_{k+1}$ from (\ref{halfinf}) we obtain with $[M,N]=[0,\infty ],$
$$b_k p_k^{-2} T_k (\CP,x) =b_k -(b_k+c_k-x)t+c_k t^2 \ge 0.$$
This yields $t_0=\frac{b_k+c_k-x}{2c_k} ,$ which is at least a $(1,k-1)-$transversal, and
\begin{equation}
\label{eqvrem}
4b_k c_k \Delta(x)= \left( x-( \sqrt{b_k}-\sqrt{c_k} \; )^2 \; \right)  
\left( ( \sqrt{b_k}+\sqrt{c_k} \; )^2-x \; \right) \ge 0.
\end{equation}
Notice that the polynomials here are normalized so that the sign of the leading coefficient 
of $p_k$ is $(-1)^k,$ and so the branches of $t(x)$ are upside-down in comparison with the previous cases.
Thus, to guarantee the intersection of $t_0$ with $B_k$
one should check 
$$\lim_{x-> \infty} t(x)=-1 > \lim_{x-> \infty} t_0(x) =- \infty ,$$
hence $t_0$ is a $(1,k)-$transversal.
On the other hand to be a $(0,k)-$transversal it should satisfy
$t_0(0)=\frac{b_k +c_k}{2 c_k} \le p_k(0)=1. $
This is the case only if $b_k \le c_k .$ This complete the proof.
\end{proof}
Note that similar bounds can be obtained for nonsymmetric polynomials orthogonal on
a finite interval.
The corresponding $\TI$ are given in \cite{szwarc}.
The restrictions we have to impose in cases (ii) and (iii) to obtain bounds on all the zeros
reflect the real situation. An easy example is provided by
ultraspherical and Laguerre polynomials with small parameters.
Roughly speaking, the reason why the obtained bounds exclude one or both extreme zeros
is that these zeros are too close to the ends of the interval of orthogonality
and the used inequalities do not have enough precision to distinguish between them.
\\
The simplest, seemingly open question, arising naturally in connection with the above theorem
is
\\
{\em Suppose that $p_k$ satisfies (\ref{recsym}), what is the maximal rate of growth
of $c_k$ such that} 
$$\lim_{k \rightarrow \infty} \frac{x_{kk}}{2 \sqrt{c_{k-1}}} =1 \; ?$$
The asymptotic (\ref{nevoz}) shows that this is true for the polynomial growth.
On the other hand it is almost certainly wrong for exponential $c_k \sim c^k, c >1.$
The following lemma supports this claim.
\begin{lemma}
\label{tt2}
Let $p_k$ be a symmetric polynomial
satisfying (\ref{recsym}) and suppose that $\frac{3}{4}c_k < c_{k-1} \le c_k .$ 
Then
$$x_{kk} < 2 \sqrt{c_{k-2}}, \; \;  k \ge 3. $$
\end{lemma}
\begin{proof}
We consider 
$$Q(t,x)=c_{k-2}c_{k-1}^2 p_k^{-2} T_{k-2}(\CP ,x)  =$$
$$
(c_{k-1}^2-(c_{k-1}-c_{k-2})x^2)t^2-x(2c_{k-2}-c_{k-1})t+c_{k-2} \ge 0.
$$
In view of Theorem \ref{firstob}, (i)  we can choose $[M,N]=[-2 \sqrt{c_{k-1}},2 \sqrt{c_{k-1}} \; ].$
Then 
$$t_0= \frac{x(2c_{k-2}-c_{k-1})}{2(c_{k-1}^2-(c_{k-1}-c_{k-2})x^2)},
$$
By the assumption $\frac{3}{4}c_k < c_{k-1} \le c_k ,$
therefore  $t_0$ is a continuous function on $[M,N]$
and thus a $(0,k)-$transversal. Finally
$$
Q(t_0,x)=\frac{(4c_{k-2}-x^2)c_{k-1}^2}{4(c_{k-1}^2-(c_{k-1}-c_{k-2})x^2)},
$$
and the result follows.
\end{proof}
Similar but more involved calculations with $T_{k-3}$ instead of $T_{k-2}$ 
yield $x_{kk} < 2 \sqrt{c_{k-3}},$ provided $\frac{5+\sqrt{5}}{8}c_k < c_{k-1} \le c_k ,$
we omit the details.
\\
Now we will show that using (\ref{marik}) which yields a fourth degree form,
one can obtain much sharper second order bounds.
We will consider the simplest case of monic Hermite polynomials $H_k$
defined by $(\ref{recsym})$ with
$c_k=k/2.$ The corresponding asymptotic for $x_{kk}$ given by (\ref{nevoz}) is
$$x_{kk} = \sqrt{2k}-2^{-1/2} 3^{-1/3} i_1 k^{-1/6} \approx \sqrt{2k}-1.65 \cdot k^{-1/6}.$$ 
Putting $u_k=H_k, \; t=H_{k-1}/H_k ,$ in (\ref{marik}), we get
$$
Q(t,x)=4 S_k(\CP, x) u^{-4}_k = 
k^2 (2k-x^2) t^4 -2 k x (1 + 4 k - 2 x^2)t^3+
$$
$$(4(k+x^2)(2k+1-x^2)-1)t^2
-4x(4k+3-2x^2)t +4(2k+2-x^2) \ge 0.
$$
Choosing the same $t_0= \frac{x}{k},$ as for the case of the quadratic form given by
$T_k $ and calculating $Q(t_0,x)$
we get that all the zeros of $p_k$ satisfy
$$
8k^2(k+1)-(6k+1)(2k+1)x^2+(6k+2)x^4-x^6 \ge 0.
$$
This equation has the only positive root $x$ being the sought bound,
$$x=\frac{(m^2-1)^2 \sqrt{m^4+4m^2+1}}{3 \sqrt{3} m^3}=\sqrt{2k}-2^{-7/6}k^{-1/6}+O(k^{-5/6}),$$
where $m=2^{-1/6} ( \sqrt{27k+2}+\sqrt{27k} \; )^{1/3}.$
Thus we get $2^{-7/6} \approx 4/9,$ instead of $1.65$. Notice that the result
can be slightly improved 
by solving the system $Q(t,x)=0,\; \frac{\partial Q(t,x)}{\partial t}=0,$
exactly. This yields 
$\sqrt{ \frac{4k-3k^{1/3}+1}{2}} \approx \sqrt{2k} -0.53 k^{-1/6},$
we omit the details.
%%%%%%%%%%%%%%%%%%%%%%%%%%%%%%%%%%%%%%%%%%%%%%%%%%%%%%%%%%%%%%%%%%

\section{Bounds from higher order Tur\'{a}n inequalities}
In this section we will consider only the symmetric case (\ref{recsym})
where we put for convenience $c_0=0.$
First, we will establish sufficient conditions for the validity of the inequality 
$$T^{(2)}_k=T_k^{(2)} (\CP,x)=3 p_k^2-4p_{k-1}p_{k+1}+p_{k-2}p_{k+2} \ge 0,$$
in terms of the recurrence (\ref{recsym}) and present the corresponding bounds on the
extreme zeros. Next, we will show how to modify $T^{(1)}_k$ to obtain second order bounds
for a vast class of nondecreasing sequences $c_k.$
In particularly, we will prove Theorems \ref{th1} and  \ref{thmnt}.
\begin{theorem}
\label{tur2}
Let $ \{ c_k \}_{k=1}^\infty , $ be a nondecreasing positive sequence such that 
\begin{equation}
\label{cond1}
c_{k-1}-3 c_k+3 c_{k+1}-c_{k+2} \ge 0.
\end{equation}
Then for $k \ge 2,$
\begin{equation}
\label{turan2}
c_{k-1}p_k^{-2} T_k^{(2)}  =
\end{equation}
$$
c_k (4c_{k-1}-x^2) t^2-
x(4 c_{k-1}-c_k+c_{k+1}-x^2)t+3 c_{k-1}+c_{k+1}-x^2 \ge 0.
$$
\end{theorem}
\begin{proof}
We have the following directly checked identity
\begin{equation}
\label{eqtur2}
T_{k+1}^{(2)} =c_{k-1}T_k^{(2)}+(c_{k+2}+3 c_k -4 c_{k-1})T_k +(c_{k-1}-3 c_k+3 c_{k+1}-c_{k+2})p_k^2 .
\end{equation}
Now the result follows by the induction on $k$ and
\\
$T_2^{(2)}  =(c_0-3c_1+3c_2-c_3)x^2+3c_1^2+c_1 c_3 > 0.$
\end{proof}
\begin{remark}
If we set $c_k=\sum_{i=1}^k \delta_i, $ then the conditions (\ref{cond1})
can be rewritten as $\delta_i \ge 0,$ $\delta_{i-1}-2 \delta_i+\delta_{i+1} <0,$
i.e. $\delta_i$ should be a nonnegative concave function of $i.$
\end{remark}
\begin{theorem}
\label{tur2oz}
Let $c_k$ satisfy the conditions (\ref{cond1}) of Theorem \ref{tur2}, then all the zeros of
$p_k$ are confined between the smallest and the largest real zeros of the equation
\begin{equation}
\label{vir1}
F(x)=x^6-2(4c_{k-1}+c_k+c_{k+1})x^4+(16 c_{k-1}^2+(c_k+c_{k+1})^2 +
\end{equation}
$$
4c_{k-1}(5c_k+2c_{k+1}))x^2-16c_k c_{k-1}(3 c_{k-1}+c_{k+1})=0.
$$
\end{theorem}
\begin{proof}
We have by (\ref{turan2}),
$$
Q(t,x)= c_{k-1}p_k^{-2} T_k^{(2)} =$$
$$3 c_{k-1}+c_{k+1}-x^2-x(4 c_{k-1}-c_k+c_{k+1}-x^2)t+
c_k (4c_{k-1}-x^2) t^2 \ge 0.
$$
Clearly, $x^2 < 4 c_{k-1}, $
and we can choose $[M,N]=[-2 \sqrt{c_{k-1}},2 \sqrt{c_{k-1}} \; ].$
Now
$$t_0=\frac{x(4 c_{k-1}-c_k+c_{k+1}-x^2)}{2c_k (4c_{k-1}-x^2)},$$
is a $(0,k)-$ transversal and any zero $x$ satisfies
$$
Q(t_0,x)=- \; \frac{F(x)}{4c_k (4c_{k-1}-x^2)} >0,
$$
yielding the required result.
\end{proof}

To show that (\ref{vir1}) indeed yields second order bounds
we again consider the monic Hermite polynomials $H_k(x).$
The conditions of Lemma \ref{tur2} are fulfilled
as $c_k=k/2.$ Solving (\ref{vir1}) we get
$$
x_{kk} < \sqrt{2k- \frac{(1+(\sqrt{k}+\sqrt{k-1} \; )^{2/3})^2}{2(\sqrt{k}+\sqrt{k-1} \; )^{2/3}}}=
\sqrt{2k-1}- 2^{-5/3}(2k-1)^{-1/6}+O(k^{-5/6}).
$$
%%%%%%%%%%%%%%%%%%%%%%%%%%%%%%%%%%%%%%%%%%%%%%%%%%%%%%%%%%%%%%%%%%%%%%%%%%%%%%%%%%%%%%%%%%%%%%%%%%%%%%%

Now we will establish a new $\TI$ which is valid for a vast class of sequences $c_k.$
\\
Given $\CP=\{ p_k \},$ define $ \CQ= \{ q_k \},$ by $ q_k(x)=p_{k+1}(x) -c_k p_{k-1}(x).$
The following identity can be checked directly.
\begin{lemma}
$$
T_{k+1}(\CQ,x)=c_k T_k (\CQ,x)+
2 c_k \mu_k T_k (\CP,x)+G,$$
where
$$ p_k^{-2} G =2c_k^2 (2 c_{k+2}-2 c_k- \mu_k )t^2 -
2x c_k (3 c_{k+2}-2c_{k+1}-c_k- \mu_k )t+ 
$$
$$x^2(2 c_{k+2}-3 c_{k+1}+c_k )+4c_{k+1}(c_{k+1}-c_k )-2 c_k \mu_k ,
$$
$\mu_k =2(c_{k+1}-c_k )+\frac{1}{2} \; (\sqrt{c_{k+1}-c_k}-\sqrt{2 c_{k+2}-3c_{k+1}+c_k} \; )^2 .$
\end{lemma}

\begin{theorem}
\label{lemmain}
Let $\{ c_k \}_{k=1}^\infty$ be a nondecreasing sequence satisfying for
$k=1,2,...,$ the following conditions
\begin{equation}
\label{condi}
2c_{k+2}-3 c_{k+1}+c_k \ge 0,
\end{equation}
\begin{equation}
\label{condii}
(c_{k+1}-c_k) ( \sqrt{c_{k+1}-c_k}+ \sqrt{2 c_{k+2}-3c_{k+1}+c_k} \; ) \ge  
\end{equation}
$$
\sqrt{c_k}\;|c_{k+2}-2c_{k+1} +c_k| .
$$
Then 
\begin{equation}
\label{vx}
p_k^{-2} T_k (\CQ,x)=c_k(4 c_k -x^2)t^2-x(2c_{k+1}+2 c_k-x^2)t+4c_{k+1}-x^2 \ge 0.
\end{equation}
\end{theorem}
\begin{proof}
As $T_k(\CP,x) \ge 0,$ by $c_{i+1} \ge c_i,$ and
$T_1(\CQ,x)=(2c_2-3c_1)x^2+4c_1^2 >0,$ by (\ref{condi}) it is left to show that $G \ge 0.$ 
For, consider
$$H=  p_k^{-2} G = 2c_k^2 (2 c_{k+2}-2 c_k- \mu_k )t^2-2x c_k (3 c_{k+2}-2c_{k+1}-c_k- \mu_k )t+$$
$$
x^2(2 c_{k+2}-3 c_{k+1}+c_k )+4c_{k+1}(c_{k+1}-c_k )-2 c_k \mu_k ,
$$
The coefficient at $t^2$ is positive, hence it is left to check that the discriminant of
this quadratic in $t$ is nonpositive,
what yields
$$c_k (c_{k+2}-2c_{k+1} +c_k)^2 -(c_{k+1}-c_k)^2 ( \sqrt{c_{k+1}-c_k}+ 
\sqrt{2 c_{k+2}-3c_{k+1}+c_k} \; )^2 \le 0,$$
and (\ref{condi}), (\ref{condii}) follow.
\end{proof}
Practically the conditions of the above theorem are much less restrictive
than these of Theorem \ref{tur2}.
Yet formally 
(\ref{cond1}) does not follow from (\ref{condi}),(\ref{condii}),
as the example $c_k=k^2+c,$ shows.
\\
The conditions (\ref{condi}) and (\ref{condii}) are rather complicated
but can be simplified by the substitution  $ c_k/c_{k-1} =1 + d_k, \; d_k >0,$ 
giving respectively
\begin{equation}
\label{condi1}
d_{k+1} \ge \frac{d_k}{2(1+d_k )} ,
\end{equation}
\begin{equation}
\label{condii1}
d_k (\sqrt{d_k}+\sqrt{2 d_k d_{k+1}+2 d_{k+1}-d_k} \; ) \ge |d_k d_{k+1}+ d_{k+1}-d_k |.
\end{equation}
More practical criteria are given in the following Lemma.
\begin{lemma}
The conditions of Theorem \ref{lemmain} hold if $d_k >0,$ and 
\begin{equation}
\label{condnew}
\frac{d_k}{2(1+d_k)} < d_{k+1} < \frac{d_k(1+2 \sqrt{d_k}+2 d_k )}{1+d_k} ,
\end{equation}
\end{lemma}
\begin{proof}
Putting in (\ref{condi1}),(\ref{condii1})  $d_{k+1} =\frac{d_k (1+2y+2y^2) }{1+d_k}, \; y \ge -1/2 ,$
we get 
$$
1+2y+2y^2 \ge \frac{1}{2}, \;\; \;  \; 4 d_k^2 (1+y)^2 (d_k-y^2) >0,
$$
and the result follows.
\end{proof}
Now we are in the position to state the main result of the paper.
\begin{theorem}
\label{thmain}
Suppose that $d_i =\frac{c_i -c_{i-1}}{c_i} \ge 0, $
satisfy (\ref{condi1}),(\ref{condii1}).
Then 
$$
x_{kk}^2<4c_k \left(1-  6^{-4/3} d_{k+1}^{2/3} \; \left(
( v+9 )^{1/3}- ( v-9  )^{1/3}
\right)^2  \right) 
$$
where $v=\sqrt{6d_{k+1}+81}.$
\end{theorem}
\begin{proof}
Let $[M,N]=[-2 \sqrt{c_k},2 \sqrt{c_k} \; ],$ and 
consider $Q(t,x)= p_k^{-2} T_k (\CQ,x)$  given by (\ref{vx}). 
Then we find
$$
t_0= \frac{x(2c_{k+1}+2 c_k-x^2)}{2c_k (4c_k-x^2)},
$$
and $t_0$ is a $(0,k)-$transversal.
Calculating $Q(t_0,x)$ we conclude that all the zeros of $p_k$ satisfy 
$$x^6-4(c_{k+1}+2c_k )x^4+4(c_{k+1}+c_k )(c_{k+1}+5c_k )x^2 -64c_k^2 c_{k+1} < 0
$$
The corresponding equation has only two real roots giving the required bounds, namely
$$
x^2=4c_k \left(1-  6^{-4/3} d_{k+1}^{2/3} \; \left(
( v+9 )^{1/3}- ( v-9  )^{1/3}
\right)^2  \right).
$$
To show that there are no other roots
we calculate its discriminant which is
$$
-2^{28} c_{k+1} c_k^4 (c_{k+1}-c_k)^8 (2c_{k+1}+25c_k)^2.
$$
As it does not change the sign, provided $c_{k+1}>c_k >0,$ the number of real zeros is the same for
any such a choice of $c_k, c_{k+1}.$
Choosing $c_k=1, c_{k+1}=2,$ we obtain the test equation
$x^6-16x^4+84x^2-128=0,$ having only two real roots.
\end{proof}

Now Theorem \ref{th1} is a direct consequence of the following claim.
\begin{lemma}
In the conditions of Theorem \ref{thmain}
\begin{equation}
\label{condsimpl1}
x_{kk}^2 < 4c_k \left( 1- \frac{d_{k+1}^{2/3}}{(2^{1/3}+d_{k+1}^{1/3} \; )^2} \right).
\end{equation}
\end{lemma}
\begin{proof}
It is enough to show that
$$
( v+9 )^{1/3}- ( v-9  )^{1/3} > \frac{6^{2/3}}{2^{1/3}+d_{k+1}^{1/3} }, 
$$
which is transformed into 
$ y > \frac{ y}{1+y-y^3}, $
by the substitution $d_{k+1}=\frac{2(1-y^3)^3}{y^3},$ $0 <y \le 1.$
\end{proof}

Finally, Theorem \ref{thmnt} follows from (\ref{condsimpl1}) with $d_{k+1}=(1+\frac{1}{k})^{2 \delta}-1.$ 
For we observe that $ \frac{d_{k+1}^{1/3}}{2^{1/3}+d_{k+1}^{1/3}}$ is an increasing
function in $d_{k+1}$ and 
and the result follows by applying the elementary inequality
$$
(1+\frac{1}{k})^{2 \delta}-1 \ge \frac{2 \delta}{k+\frac{1}{2}}, \; \; \delta \ge 0. 
$$
Moreover, $k+\frac{1}{2}$ may be replaced by $k$ for $\delta \ge \frac{1}{2}.$ 
%%%%%%%%%%%%%%%%%%%%%%%%%%%%%%%%%%%%%%%%%%%%%%%%%%%%%%%%%%%%%%%%%%%%%%%%%%%%%%%%%%%%%%%%%%%%%%%%%%%%%%%%%%%%


\begin{thebibliography}{8}

\bibitem{assche00}
D.W.\,Lee, W.\,Van\,Assche, {\em Asymptotic of orthogonal polynomials by three term
recurrence relation}, preprint.
%\bibitem{assche}
%W.\,Van\,Assche and A.\,B.\,J.\,Kuijlaars,
%{\em The asymptotic distribution of orthogonal
%polynomials with varying recurrence coefficients},
%J. Approx. Theory, 99, (1999), 167-197.
\bibitem{craven}
T.\,Craven, G.\,Csordas, {\em Iterated Tur\'{a}n and Laguerre inequalities},
J. Ineq. in Pure and Appl. Math. 3, (2003). 
\bibitem{craven1}
T.\,Craven, G.\,Csordas, {\em Composition theorems, multiplier sequences and
complex zero decreasing sequences}, preprint.
\bibitem{dimit}
D.K.\,Dimitrov, {\em Higher order Tur\'{a}n inequalities}, Proc. Amer. Math. Soc. 126 (1998), 2033-2037.
\bibitem{freud}
G.\,Freud, {\em On the greatest zero of an orthogonal polynomial}, J.Approx. Theory 46 (1986), 16-24.
\bibitem{ismail} M.E.H\,Ismail and X.\,Li, Bounds on the extreme zeros of orthogonal
polynomials, {\em Proc. Amer. Math. Soc.} 115, 1992, 131-140.
\bibitem{k3}
I.\,Krasikov, {\em Nonnegative quadratic forms and bounds on
orthogonal polynomials}, J. Approx. Theory 111, (2001), 31-49.
\bibitem{klag}
I.\,Krasikov, {\em Bounds for zeros of the Laguerre polynomials}, 
J. Approx. Theory, 121, (2003) 287-291.
\bibitem{kest}
I.\,Krasikov, {\em On zeros of polynomials and allied functions satisfying
second order differential equation},
East J. Approx., 9 (2003) 51-65.
\bibitem{kd}
I.\,Krasikov, {\em Discrete analogues of the Laguerre inequality},
Analysis and Applications, 1 (2003), 189-198.
\bibitem{kch}
I.\,Krasikov, {\em Bounds for zeros of the Charlier polynomials},
Methods and Applications of Analysis, to appear.
\bibitem{kf}
I.\,Krasikov, {On extreme zeros of classical orthogonal polynomials}, submitted.
\bibitem{leven} 
V.\,I.\,Levenstein,{\em  Universal bounds on codes and designs}, 
In: Handbook of Coding Theory, Vol.1, North-Holland, 1998, 499-648.
\bibitem{marik}
J.\,Ma\v{r}\'{\i}k, {\em On polynomials with all real zeros}, \v{C}asopis P\v{e}st. Mat. 89 (1964), 5-9.
\bibitem{mtn}
A.\,M\'{a}t\'{e}, P.\,Nevai, V.\,Totik, {\em Asymptotic of 
the zeros of orthogonal polynomials associated with infinite
intervals}, J. London. Math. Soc. 33, (1986), 303-310.
\bibitem{nicu}
C.P.\,Niculescu, {\em A new look at Newton's inequalities},
J. Ineq. in Pure and Appl. Math. 1, (2000).
%\bibitem{patrick1}
%M.\,L.\,Patrick, Some inequalities concerning Jacobi Polynomials,
%{\em SIAM J. Math. Anal.}, 2, 1971, 213-220.
\bibitem{patrick2}
M.\,L.\,Patrick,  Extension of inequalities of the Laguerre and Tur\'{a}n type,
{\em Pacific J. Math.} 44, (1973), 675-682.
%\bibitem{szego} G.\,Szeg\"{o},
%{\em Orthogonal Polynomials},
%Amer. Math. Soc. Colloq. Publ., v.23, Providence, RI, 1975.
\bibitem{szwarc}
R.\,Szwarc, {\em Positivity of Tur\'{a}n determinants for orthogonal polynomials},
in {\em Harmonic Analysis and Hypergroups}, (K.A.\,Ross et al., ed.) Delhi 1995, Birkhauser,
Boston-Basel-Berlin 1997, 165-182.
\bibitem{turan}
P.\,Tur\'{a}n, {\em On the zeros of the polynomials of Legendre},
\v{C}asopis P\v{e}st. Mat. 75, (1950), 113-122
\end{thebibliography}
\end{document}